\documentclass[11pt,openbib]{article}
\usepackage{amsmath}
\usepackage{amssymb}
\usepackage{amsfonts}
\usepackage[T1]{fontenc}
\usepackage{tikz-cd}

\newcommand{\R}{\mathbb{R}}

\newcommand{\Gl}{\mathrm{Gl}}

\newcommand{\dee}{\mathrm{d}}
\newcommand{\comp}{\, \raisebox{2pt}{$\scriptstyle\circ \, $}}

\newcommand{\spann}{\mathop{\rm span}\nolimits}

\newcommand{\dbydt}{\mbox{${\displaystyle \frac{\dee }{\dee t}}
\rule[-10pt]{.5pt}{25pt} \raisebox{-10pt}{$\, {\scriptstyle t=0}$}$}}

\newcommand{\setrule}{\, \rule[-4pt]{.5pt}{13pt}\, }

\allowdisplaybreaks

\begin{document}
\thispagestyle{empty}
\begin{center}
{\Large \textbf{Surjective submersions of subcartesian spaces}} \\ 
\mbox{}\vspace{-.2in} \\
\mbox{}\\
Richard Cushman
\mbox{}\\ 
\date{} 
\end{center}
\addtocounter{footnote}{1}
\footnotetext{printed: \today }
\begin{abstract}
We define the notion of a submersion of subcartesian differential spaces and 
prove some of its properties, which are analogous to those of a submersion in 
the category of smooth manifolds and smooth mappings.
\end{abstract}\bigskip

\noindent{\Large \textbf{Introduction}} \medskip 

The aim of this paper is to define the notion of a surjective submersion of subcartesian 
differential spaces. We show that a surjective submersion of subcartesian spaces is an open mapping. \medskip 

\section{Aronszajn and Sikorski subcartesian spaces}

We begin by reviewing some results on Aronszajn and Sikorski subcartesian spaces, 
which are proved in \cite{cushman-sniatycki21}. \medskip 

A differential space $(S, C^{\infty}(S))$ is \emph{subcartesian} if its differential space 
topology is Hausdorff and for every point $x \in S$ there is an open set 
$U_x$ of $x$ in $S$ and a one to one mapping ${\varphi }_x: U_x \subseteq S \rightarrow 
V_x = {\varphi }_x(U_x) \subseteq {\R }^{n_{{\varphi }_x}}$ which is a diffeomorphism 
of the differential space $(U_x, C^{\infty}(U_x))$ onto the differential subspace $(V_x, C^{\infty}(V_x))$ 
of $( {\R }^{n_{{\varphi }_x}}, C^{\infty}({\R }^{n_{{\varphi }_x}}))$. \medskip 

A \emph{subcartesian space} $(S, \mathfrak{A})$ \emph{in the sense of Aronszajn} is a Hausdorff 
topological space $S$ endowed with an atlas $\mathfrak{A} = 
\{ \varphi : U_{\varphi } \rightarrow V_{\varphi } \}$, 
where $\varphi $ is a homomorphism of an open set $U_{\varphi }$ of $S$ onto a subset 
$V_{\varphi }$ of ${\R }^{n_{\varphi}}$ which satisfies \medskip 

1. \parbox[t]{4.25in}{The domains $\{ U_{\varphi } \setrule \, \varphi \in \mathfrak{A} \}$ form an 
open covering of $S$.} \smallskip 

2. \parbox[t]{4.25in}{For every $\varphi $, $\psi \in \mathfrak{A}$ and every $x \in U_{\varphi }\cap 
U_{\psi}$ there is a $C^{\infty}$ mapping ${\Phi }_{\varphi }$ in an open neighborhood $W_{\varphi (x)}$ of $\varphi (x)$ in ${\R }^{n_{\varphi }}$, which extends the mapping $\psi \comp {\varphi }^{-1}: 
\varphi (U_{\varphi } \cap U_{\psi }) \rightarrow \psi (U_{\varphi } \cap U_{\psi })$ and a $C^{\infty}$ 
mapping ${\Phi }_{\psi }$ in an open neighborhood $W_{\psi (x)}$ of $\psi (x)$ in ${\R }^{n_{\psi}}$ which extends the mapping $\varphi \comp {\psi }^{-1}: \psi (U_{\varphi } \cap U_{\psi }) \rightarrow \varphi (U_{\varphi } \cap U_{\psi })$.} 
\bigskip 

Let $(S_1, {\mathfrak{A}}_1)$ and $(S_2, {\mathfrak{A}}_2)$ be subcartesian spaces in the sense of 
Aronszajn. A mapping $\chi : S_1 \rightarrow S_2$ is \emph{smooth} if for every $y \in S_1$ there 
are ${\varphi }_1 \in {\mathfrak{A}}_1$ and ${\varphi }_2 \in {\mathfrak{A}}_2$ such that $y \in U_{{\varphi }_1}$ 
and for every $x \in U_{{\varphi }_1}$ we have $\chi (x) \in U_{{\varphi }_2}$ and the mapping 
${\varphi }_2 \comp \chi \comp {\varphi }^{-1}_1: V_{{\varphi }_1} \rightarrow V_{{\varphi }_2}$ 
extends to a $C^{\infty}$ mapping ${\Phi }_{{\varphi }_1, {\varphi }_2}$ of an open 
neighborhood $W_{{\varphi}_1(x)}$ of ${\varphi }_1(x)$ in 
${\R }^{n_{ {\varphi }_1}}$ to an open neighborhood $W_{{\varphi }_2(\chi (x))}$ of ${\varphi }_2(\chi (x))$ in 
${\R }^{n_{{\varphi }_2}}$. \medskip 

\noindent \textbf{Proposition 1.} Let $(S, C^{\infty}(S))$ be a subcartesian differential space with atlas 
$\mathfrak{A} = {\{ {\varphi }_x: U_x \rightarrow V_x \}}_{x\in S}$. Then $(S, \mathfrak{A})$ is 
a subcartesian space in the sense of Aronszajn. \medskip 

The atlas $\mathfrak{A}$ in proposition 1 is called the \emph{atlas associated to 
the differential structure} $C^{\infty}(S)$. \medskip 

Let $(S, \mathfrak{A})$ be a subcartesian space in the sense of Aronszajn. A function 
$f: (S, \mathfrak{A}) \rightarrow \R $ lies in the family $\mathcal{F}$ if for every 
$y \in S$ there is $\varphi \in \mathfrak{A}$ such that $y \in U_{\varphi }$ and if $x \in U_{\varphi }$ then 
there is an open neighborhood $W_{\varphi (x)}$ of $\varphi (x) \in {\R }^{n_{\varphi }}$ 
and a function $F_{\varphi (x)} \in C^{\infty}({\R }^{n_{\varphi }})$ satisfying 
$f = F_{\varphi (x)} \comp \varphi $ on ${\varphi }^{-1}( V_{\varphi } \cap W_{\varphi (x)})$. 
There is a differential structure $C^{\infty}(\mathfrak{A})$ on $S$ generated by the family $\mathcal{F}$. This 
differential structure is \emph{determined by the atlas} $\mathfrak{A}$ \emph{on} $S$. The spaces 
$(S, \mathfrak{A})$ and $(S, C^{\infty}(\mathfrak{A}))$ have the same topology. \medskip %

\noindent \textbf{Proposition 2.} Let $(S_1, {\mathfrak{A}}_1)$ and $(S_2, {\mathfrak{A}}_2)$ be 
Aronszajn subcartesian spaces and let $(S_1, C^{\infty}({\mathfrak{A}}_1))$ and $(S_2, 
C^{\infty}({\mathfrak{A}}_2))$ be subcartesian differential spaces whose differential structure is determined by the corresponding atlas. If $\chi : (S_1, {\mathfrak{A}}_1) \rightarrow (S_2, {\mathfrak{A}}_2)$ is a smooth map 
of Aronszajn subcartesian spaces, then $\chi : (S_1, C^{\infty}({\mathfrak{A}}_1)) \rightarrow 
(S_2, C^{\infty}({\mathfrak{A}}_2))$ is a smooth map of subcartesian differential spaces. 

\section{Submersion}

In this section we define the notion of a submersion of subcartesian differential spaces. \medskip 

Let $\chi : (S_1, C^{\infty}(S_1)) \rightarrow (S_2, C^{\infty}(S_2))$ be a surjective smooth map 
of subcartesian differential spaces. Let ${\mathfrak{A}}_1$ and ${\mathfrak{A}}_2$ be 
atlases associated to the differential structures $C^{\infty}(S_1)$ and $C^{\infty}(S_2)$, 
respectively. We say that 
$\chi : (S_1, C^{\infty}(S_1)) \rightarrow (S_2, C^{\infty}(S_2))$ is a \emph{submersion} if 
for every $y \in S_1$ there is a ${\varphi }_1 \in {\mathfrak{A}}_1$ and a ${\varphi }_2 \in {\mathfrak{A}}_2$ 
such that $y \in U_{{\varphi }_1}$ and for every $x \in U_{{\varphi }_1}$ we have $\chi (x) \in U_{{\varphi }_2}$ and the mapping ${\varphi }_2 \comp \chi \comp {\varphi}^{-1}_1: V_{{\varphi }_1} \rightarrow V_{{\varphi }_2}$ extends 
to a $C^{\infty}$ submersion ${\Phi }_{{\varphi }_1, {\varphi }_2}$ of an open neighborhood 
$W_{{\varphi }_1(x)}$ of ${\varphi }_1(x)$ in 
${\R }^{n_{{\varphi }_1}}$ onto an open neighborhood $W_{{\varphi }_2(\chi (x))}$ of ${\varphi }_2(\chi (x))$. 
In other words, the map $\chi : (S_1, {\mathfrak{A}}_1) 
\rightarrow (S_2, {\mathfrak{A}}_2)$ is a submersion of Aronszajn subcartesian spaces. 
Conversely, if  $\chi : (S_1, {\mathfrak{A}}_1) \rightarrow (S_2, {\mathfrak{A}}_2)$ is a submersion 
of Aronszajn subcartesian spaces, then $\chi : (S_1, C^{\infty}({\mathfrak{A}}_1)) \rightarrow 
(S_2, C^{\infty}({\mathfrak{A}}_2))$ is a submersion of subcartesian differential spaces.  \medskip 

\noindent \textbf{Theorem 3.} Suppose that $\chi : (S_1, C^{\infty}(S_1)) \rightarrow (S_2, C^{\infty}(S_2))$ is 
a surjective submersion of subcartesian differential spaces. Then $\chi $ is an open mapping, 
that is, if $U$ is an open subset of $S_1$, then $\chi (U)$ is an open subset of $S_2$. \medskip 

\noindent \textbf{Proof.} Let $U$ be an open subset of $S_1$. Since 
$\chi $ is a surjective submersion for each $x \in U$ there are diffeomorphisms ${\varphi }_x \in 
{\mathfrak{A}}_1$ and ${\varphi }_{\chi (x)} \in {\mathfrak{A}}_2$ with the following properties. \medskip 

\indent 1. \parbox[t]{4.25in}{The atlases ${\mathfrak{A}}_1$ and ${\mathfrak{A}}_2$ are determined by the 
differential structures $C^{\infty}(S_1)$ and $C^{\infty}(S_2)$, respectively.} \smallskip 

\indent 2. \parbox[t]{4.25in}{The diffeomorphisms ${\varphi }_x: U_x \rightarrow V_x \subseteq 
{\R }^{ n_{ {\varphi }_x} }$ and ${\varphi }_{\chi (x) }: U_{\chi (x)} \rightarrow V_{\chi (x)} \subseteq 
{\R }^{ n_{ {\varphi }_{\chi (x)} } }$ lie in the atlases ${\mathfrak{A}}_1$ and ${\mathfrak{A}}_2$, 
respectively.} \smallskip  

\indent 3. \parbox[t]{4.25in}{For every $\xi \in V_x$ the map 
\begin{displaymath}
{\varphi }_{\chi (x)} \comp \chi \comp {\varphi }^{-1}_x: V_x \subseteq 
{\R }^{n_{{\varphi}_x}} \rightarrow V_{\chi (x)} \subseteq {\R }^{n_{{\varphi }_{\chi (x)} }}
\end{displaymath}
extends to a submersion ${\Phi }_{\xi }$ of an open subset $V^{\xi}$ of ${\varphi }_x(\xi )$ in 
${\R }^{n_{{\varphi }_x}}$ onto an open subset $V^{\eta }$ of ${\Phi }_{\xi }({\varphi}_x(\xi )) = 
{\varphi }_{\chi (x)}(\chi (\xi ))$ in ${\R }^{n_{{\varphi }_{\chi (x)}}}.$} \medskip 

\noindent So $({\varphi }_{\chi (x)} \comp \chi \comp {\varphi }^{-1}_x)(V^{\xi } \cap V_x) = V^{\eta } \cap V_{\chi (x)}$. In other words,  
\begin{equation}
\chi \big( {\varphi }^{-1}_x(V^{\xi } \cap V_x) \big) = {\varphi }^{-1}_{\chi (x)}( V^{\eta } \cap V_{\chi (x)}) 
\label{eq-one}
\end{equation} 
is an open subset of $\chi (U_x) \subseteq U_{\chi (x)}$ containing $\chi (\xi )$. Letting $\xi $ range over $V_x$, equation (\ref{eq-one}) implies 
\begin{displaymath}
\chi (U_x) = \chi \big( {\varphi }^{-1}_x(V_x) \big) = 
{\varphi }^{-1}_{\chi (x)}(V_{\chi (x)}) = U_{\chi (x)}, 
\end{displaymath}
which is a open subset of $\chi (U)$ containing $\chi (x)$. Consequently, 
$\chi (U)$ is an open subset of $S_2$. \hfill $\square $ \medskip 

\section{An infinitesimal criterion}

In this section we give an infinitesimal criterion for a smooth surjective mapping of 
subcartesian differential spaces to be a submersion. \medskip

Let $(S, C^{\infty}(S))$ be a subcartesian differential space. A 
\emph{vector field} on $S$ is a derivation $X$ of $C^{\infty}(S)$ such that for every 
$x \in S$ there is an open neighborhood $U$ of $x$ in $S$ and an $\varepsilon >0$ such 
that for every $t \in (-\varepsilon , \varepsilon )$ the map 
${\varphi }_t: U \rightarrow V_t = {\varphi }_t(U) \subseteq S$, where $t \mapsto {\varphi }_t(x)$ 
is an integral curve of $X$, is defined and is a diffeomorphism onto the open subset $V_t$ of $S$. \medskip 

\noindent \textbf{Lemma 4.} The set $\mathfrak{X}(S)$ of all vector 
fields on a subcartesian space $S$ is a real vector subspace of 
$\mathrm{Der}\, C^{\infty}(S)$, the set of derivations of $C^{\infty}(S)$. \medskip 

\noindent \textbf{Proof.}
 Let $X$, $Y \in \mathfrak{X}(S)$ and $\alpha $, $\beta \in \R $. Let $\mathcal{F}$ be the smallest 
 Lie subalgebra of $(\mathfrak{X}(S), [\, \, , \, \, ])$ containing the vector fields $X$ and 
 $Y$. Let $x \in S$. Then the orbit ${\mathrm{O}}_x$ of $\mathcal{F}$ through $x$ is a smooth submanifold of $S$. Hence $X_{| {\mathrm{O}}_x}$ and 
$Y_{|{\mathrm{O}}_x}(x)$ are vector fields on the smooth manifold ${\mathrm{O}}_x$, 
see \cite{sniatycki13}. Let ${\gamma }_x:I_x 
\subseteq \R \rightarrow {\mathrm{O}}_x: t \mapsto 
{\varphi }^{\alpha X_{|{\mathrm{O}}_x} + \beta Y_{|{\mathrm{O}}_x}}_t(x)$ be an integral curve of the vector field 
$(\alpha X + \beta Y)_{|{\mathrm{O}}_x}$ on 
${\mathrm{O}}_x$ starting at $x$. Then ${\gamma }_x$ is an integral 
curve of the derivation $\alpha X + \beta Y$ of $C^{\infty}(S)$ starting at 
$x$. Hence $I_x$ contains an open interval 
about $0$ where ${\gamma }_x$ is defined. Suppose that $I_x$ is the maximal 
interval where ${\gamma }_x$ is defined. If $I_x = (-\infty , \infty)$ for every $x \in S$, then 
$\alpha X + \beta Y$ is a vector field on $S$. Suppose that there is an $x \in S$ such that 
${\overline{I}}_x$ has finite right end point $t_x$. Let $y = 
{\varphi }^{\alpha X + \beta Y}_{t_x}(x)$. The vector field 
$(\alpha X + \beta Y)|_{{\mathrm{O}}_y}$ is 
defined at $y$ with $t \mapsto {\varphi }^{\alpha X + \beta Y}_t(y)$ being an integral curve starting at $y$, which agrees with ${\gamma }_x$, when $t \in I_x$ and $t < 0$. This integral curve is defined on an open interval containing $0$, which contradicts the definition of $t_x$. Hence 
the right end point $t_x$ of ${\overline{I}}_x$ does not lie in $I_x$. A similar 
argument shows that the left end point of ${\overline{I}}_x$ does not lie 
in $I_x$. Hence $I_x$ is an open interval, which implies that the derivation 
$\alpha X + \beta Y$ on $C^{\infty}(S)$ is a vector field on $S$, since 
$S$ is a subcartesian space. \hfill $\square $ \medskip 

For each $x \in S$ the \emph{strong tangent space} ${\mathcal{T}}_xS$ of $S$ is 
the vector space subspace ${\spann }_{\R} \{ X(x) \in {\mathrm{Der}}_x C^{\infty}(S) \setrule $ 
$X \in \mathfrak{X}(S) \}$ of $T_xS$. Recall that the \emph{tangent space} $T_xS$ to 
$S$ at $x\in S$ is the set of derivations of $C^{\infty}(S)$. Let $TS = \coprod_{x\in S} T_x$ be the 
\emph{tangent bundle} of $S$. The \emph{tangent} $T\varphi : TR \rightarrow TS$ of a smooth 
map $\varphi :R \rightarrow S$ of differential spaces sends the derivation $\nu $ of $C^{\infty}(R)$ 
to the derivation $T\varphi (\nu )$ of $C^{\infty}(S)$ given by $\nu ({\varphi }^{\ast }\, \, )$. 
Let $\chi : R \rightarrow S$ be a smooth mapping of subcartesian differential spaces. Suppose that for each vector field $X \in \mathfrak{X}(R)$ and each $x \in R$ there is a vector field $Y \in \mathfrak{X}(S)$, which is 
$\chi $ related to $X$ at $\chi (x)$, that is, 
\begin{displaymath}
\dbydt \hspace{-9pt}(\chi \comp {\varphi }^X_t)(x) = \dbydt \hspace{-9pt}{\varphi }^Y_t \big( \chi (x) \big) . 
\end{displaymath}
Then 
\begin{displaymath}
{\mathcal{T}}_x\chi : {\mathcal{T}}_xR \rightarrow {\mathcal{T}}_{\chi (x)}S: 
X(x) \mapsto Y\big( \chi (x) \big) 
\end{displaymath} 
is the \emph{strong tangent map} of $\chi $ at $x$. In other words, the map $\chi : R \rightarrow S$ has 
a strong tangent map $\mathcal{T}\chi : \mathcal{T}R \rightarrow \mathcal{T}S$ if and only if 
${\iota }_{\mathcal{T}S} \comp \mathcal{T}\chi = T \chi \comp {\iota }_{\mathcal{T}R}$. Here 
${\iota }_{\mathcal{T}S}: \mathcal{T}S \rightarrow TS$ and ${\iota }_{\mathcal{TR}}:\mathcal{T}R \rightarrow 
TR$ are inclusion mappings.  \medskip 

Let $f \in C^{\infty}(S)$. The function $f$ is \emph{strongly smooth} if the smooth mapping 
$f: (S, C^{\infty}(S)) \rightarrow (\R , C^{\infty}(\R ) )$ of subcartesian differential 
spaces has a strong tangent.  \medskip 

\noindent \textbf{Lemma 5.} Every $f \in C^{\infty}(S)$ is strongly smooth.\footnote{This lemma was 
pointed out to me by J. \'{S}niatycki \cite{sniatycki23}.} \medskip  

\noindent \textbf{Proof.} Suppose that $x \in S$ and $X(x) \in {\mathcal{T}}_xS$, where 
$X \in \mathfrak{X}(S)$. Since $f$ is smooth, $(T_xf)X(x)$ is the derivation 
$Y \big( f(x) \big) = X(f)(x)$ at $f(x)$ of $C^{\infty}(\R )$. But $(\R , C^{\infty}(\R ))$ is a smooth manifold. 
So $Y$ is a vector field on $\R $. Since $(T_xf)X(x) = Y\big( f(x) \big) $ for every $x \in S$, the 
vector fields $X$ and $Y$ are $f$-related. Hence $f$ has a strong tangent map  
\begin{displaymath}
\mathcal{T}f: \mathcal{T}S \rightarrow \mathcal{T}\R = T\R: X(x) \mapsto Y\big( f(x) \big) .
\end{displaymath} 
So $f$ is a strongly smooth function on $S$. \hfill $\square $ \medskip 

\noindent \textbf{Proposition 6.} Suppose that the differential space 
$(S, C^{\infty}(S))$ is a differential subspace of $({\R }^n, C^{\infty}({\R }^n))$. Then the inclusion mapping 
$\iota : S \rightarrow {\R }^n$ has a strong tangent map, that is, ${\mathcal{T}}_xS \subseteq T_x{\R }^n$ 
for every $x \in S$. \medskip 

\noindent \textbf{Proof.} Let $X(x) \in {\mathcal{T}}_xS$ where $x\in S$ and $X \in \mathfrak{X}(S)$. 
The vector field $X$ has a maximal integral curve 
${\gamma }_x: I \subseteq \R \rightarrow S$ starting at $x$, which is defined on a maximal open interval 
$I$ containing $0$. So ${\Gamma}_{\iota (x)} = \iota \comp {\gamma }_x: I \subseteq \R \rightarrow {\R }^n$ 
is a maximal integral curve of the vector field $Y$ on ${\R }^n$, where 
\begin{equation}
\dbydt (\iota \comp {\varphi }^X_t)(x) = Y\big( \iota (x) \big) . 
\label{eq-two}
\end{equation}
Thus ${\mathcal{T}}_x\iota X(x) = Y\big( \iota (x) \big) \in {\mathcal{T}}_{\iota (x)}{\R }^n = T_x{\R }^n$, 
which implies that the inclusion mapping $\iota : S \rightarrow {\R }^n$ has a strong tangent map. Since 
${\mathcal{T}}_x \iota : {\mathcal{T}}_xS \rightarrow T_x {\R }^n$ is the inclusion mapping and $X(x) \in T_x{\R }^n$, we obtain ${\mathcal{T}}_xS \subseteq T_x{\R }^n$. \hfill $\square $ \medskip 

Suppose that $(S, C^{\infty}(S))$ is a subcartesian space with 
$(S, \mathfrak{A})$ the corresponding Aronszajn subcartesian space. 
$\mathfrak{A}$ is the atlas ${\{ {\varphi }_x: U_x \rightarrow V_x \} }_{x\in S}$  
associated to the differential structure $C^{\infty}(S)$. Suppose that 
every ${\varphi }_x \in \mathcal{A}$ is a smooth diffeomorphism of 
$(U_x, C^{\infty}(U_x))$ onto $(V_x, C^{\infty}(V_x) )$ with a strong tangent map. 
Then $\mathfrak{A}$ is a \emph{strong atlas} of $S$, which we denote by 
${\mathfrak{A}}_{\mathrm{stg}}(S)$. We say that 
the subcartesian space $S$ with differential structure 
$C^{\infty}({\mathfrak{A}}_{\mathrm{stg}}(S))$ \linebreak 
determined by the atlas ${\mathfrak{A}}_{\mathrm{stg}}(S)$ is 
a \emph{strong} subcartesian space, which we \linebreak 
denote by $(S, C^{\infty}({\mathfrak{A}}_{\mathrm{stg}}(S)))$.  \medskip 

\noindent \textbf{Proposition 7.} Suppose that $\chi : \big( R, C^{\infty}({\mathfrak{A}}_{\mathrm{stg}}(R)) \big) 
\rightarrow \big(S, C^{\infty}({\mathfrak{A}}_{\mathrm{stg}}(S)) \big) $ is a smooth mapping of strong subcartesian spaces with a strong tangent map. If at each point $x \in R$ the strong tangent 
${\mathcal{T}}_x\chi :{\mathcal{T}}_xR \rightarrow {\mathcal{T}}_{\chi (x)}S$ of $\chi $ is surjective, then $\chi $ is a submersion. \medskip   

\noindent \textbf{Proof.} For every $x \in R$ there are diffeomorphisms 
${\varphi}_x: U_x \rightarrow V_x$ in ${\mathfrak{A}}_{\mathrm{stg}}(R)$ and 
${\varphi }_{\chi (x)}: U_{\chi (x)} \rightarrow V_{\chi (x)}$ in ${\mathfrak{A}}_{\mathrm{stg}}(S)$ 
with $\chi (U_x) \subseteq U_{\chi (x)}$ such that 
\begin{displaymath}
{\Psi }_{x, \chi (x)} = {\varphi }_{\chi (x)} \comp \chi \comp {\varphi}^{-1}_x: 
V_x \subseteq {\R}^{n_{{\varphi }_x}} \rightarrow V_{\chi (x)} \subseteq {\R }^{n_{{\varphi }_{\chi (x)}}} 
\end{displaymath}
is a smooth map of $(V_x, C^{\infty}(V_x))$ into $(V_{\chi (x)}, C^{\infty}(V_{\chi (x)}))$, which are 
differential subspaces of $({\R}^{n_{{\varphi }_x}}, C^{\infty}({\R}^{n_{{\varphi }_x}}))$ and 
$({\R }^{n_{{\varphi }_{\chi (x)}}}, C^{\infty}({\R }^{n_{{\varphi }_{\chi (x)}}}))$, respectively. Since 
the maps ${\varphi }_x$, $\chi $, and ${\varphi }_{\chi (x)}$ are smooth maps with strong tangent maps, the map 
${\Psi }_{x, \chi (x)}$ has a strong tangent map 
\begin{displaymath}
{\mathcal{T}}_y{\Psi }_{x, \chi (x)} = T_y{\Psi }_{x, \chi (x)}: 
{\mathcal{T}}_yV_x = T_y{\R}^{n_{{\varphi }_x}} \rightarrow {\mathcal{T}}_{\chi (y)}V_{\chi (x)} = 
T_{\chi (y)}{\R }^{n_{{\varphi }_{\chi (x)}}} 
\end{displaymath}
for every $y \in V_x$. By hypothesis $T_y {\Psi }_{x, \chi (x)}$ is surjective and 
hence ${\Psi }_{x, \chi (x)}$ is a submersion. Thus the mapping $\chi $ is a submersion. \hfill $\square $ \medskip 

We give a criterion for the fibration determined by a submersion to be locally trivial. In the category 
of smooth manifolds and smooth maps we have \medskip 

\noindent \textbf{Lemma 8.} Suppose that $F: M \rightarrow N$ is a surjective submersion of 
smooth manifolds with the following properties. \medskip 

1. \parbox[t]{4.25in}{For every $n' \in N$ there is an open neighborhood $V$ of $n'$ in $N$ 
such that for every $n \in V$ there is a diffeomorphism  
${\varphi }_n: F^{-1}(n) \rightarrow F^{-1}(n')$.} \smallskip 

2. \parbox[t]{4.25in}{The map $n \mapsto {\varphi }_n$ is smooth. In other words, the map 
\begin{equation}
\begin{array}{l}
\Phi : U = F^{-1}(V) =\coprod_{v\in V}F^{-1}(v) \rightarrow V \times F^{-1}(n'): \\
\rule{0pt}{13pt}\hspace{.75in} u = (v,y) \mapsto \big( v, {\varphi }_v(y) \big) 
\end{array}
\label{eq-onevnw}
\end{equation}
is smooth.} \bigskip 

\noindent Then the fibration defined by $F$ is locally trivial, that is, 
for the open neighborhood $V$ given in the first point $F^{-1}(V)$ is diffeomorphic to $V \times F^{-1}(n')$. \medskip 

\vspace{-.15in}\noindent \textbf{Proof.} The following argument shows that the smooth mapping 
$\Phi $ (\ref{eq-onevnw}) is a local diffeomorphism. The tangent map at $u \in U$ to $\Phi $ is 
\begin{displaymath}
\begin{array}{l}
T_u\Phi : T_u U = T_vV \times T_y\big( F^{-1}(v) \big) \rightarrow 
T_vV \times T_{{\varphi }_v(y)}\big( F^{-1}(n') \big):  \\ 
\rule{0pt}{16pt}\hspace{1in} \mbox{\footnotesize $\begin{pmatrix} \xi \\ \eta \end{pmatrix}$} 
\mapsto \mbox{\footnotesize $\begin{pmatrix} I & D_1{\varphi }_v(y) \\ 0 & D_2{\varphi }_v \end{pmatrix} \,  
\begin{pmatrix} \xi \\ \eta \end{pmatrix}. $}
\end{array}
\end{displaymath}
For fixed $v \in V$, the mapping ${\varphi }_v$ is a diffeomorphism by hypothesis. Hence the 
partial derivative $D_2{\varphi }_v$ in the $T_yF^{-1}(v)$ direction is invertible. Thus 
$T_u\Phi $ is invertible. So $\Phi $ (\ref{eq-onevnw}) is a local diffeomorphism. In 
fact $\Phi $ is a diffeomorphism because 1) $\Phi $ is surjective, since 
${\pi }_1 \comp \Phi  = F|U$ and $F$ is surjective. Here ${\pi }_1$ is projection on the first 
factor of $V \times F^{-1}(n')$; 2) the map $\Phi $ is injective. For if $\Phi (u) = \Phi (u')$ for 
some $u = (v,y)$, $u' = (v',y') \in U$, then $\big( v, {\varphi }_v(y) \big) = \big( v', {\varphi }_{v'}(y') \big) $, 
which gives $v = v'$. Thus ${\varphi }_v(y) = {\varphi }_{v'}(y') = {\varphi }_v(y')$. But 
${\varphi }_v$ is a diffeomorphism by hypothesis. So $y = y'$, which shows that $u =u'$. 
Thus the map $\Phi $ is injective and hence is a diffeomorphism. \hfill $\square $ \medskip 

We now prove an analogue of lemma 8 in the category of strong subcartesian spaces and  
smooth mappings with a strong tangent. \medskip 

\noindent \textbf{Theorem 9.} Suppose that $\chi : (R, C^{\infty}({\mathfrak{A}}_{\mathrm{stg}}(R))) \rightarrow 
(S, C^{\infty}({\mathfrak{A}}_{\mathrm{stg}}(S)))$ 
is a smooth surjective submersion of subcartesian spaces with a strong tangent map having the following properties. \medskip 

1. \parbox[t]{4.25in}{For every $s' \in S$ there is an open neighborhood $V$ of $s' \in S$ 
such that for every $s \in V$ there is a diffeomorphism 
\begin{displaymath}
     {\varphi }_s: \big( {\chi }^{-1}(s), C^{\infty}({\mathfrak{A}}_{\mathrm{stg}}({\chi }^{-1}(s))) \big) 
     \rightarrow \big( {\chi }^{-1}(s'), C^{\infty}({\mathfrak{A}}_{\mathrm{stg}}({\chi }^{-1}(s'))) \big). 
     \end{displaymath} } \newline 
\mbox{}\indent 2. \parbox[t]{4.25in}{The smooth map $s \mapsto {\varphi }_s$ has a strong tangent. In other 
words, the smooth map 
\begin{equation}
\begin{array}{l}
\Psi : U= {\chi }^{-1}(V) =\coprod_{v\in V}{\chi }^{-1}(v) \rightarrow V \times {\chi }^{-1}(s'): \\
\rule{0pt}{13pt}\hspace{.75in} u = (v,y) \mapsto \big( v, {\varphi }_v(y) \big) 
\end{array}
\label{eq-twovnw}
\end{equation}
of the strong subcartesian differential spaces $\big( U, C^{\infty}({\mathfrak{A}}_{\mathrm{stg}}(U)) \big) $ and 
$\big( V \times {\chi }^{-1}(s'), C^{\infty}({\mathfrak{A}}_{\mathrm{stg}}(V \times {\chi }^{-1}(s'))) \big)$ has a  
strong tangent.} \bigskip 

\noindent Then the fibration defined by $\chi $ is locally trivial, that is, 
for the open neighborhood $V$ given in the first point  
$\big( {\chi }^{-1}(V), C^{\infty}({\mathfrak{A}}_{\mathrm{stg}}({\chi }^{-1}(V))) \big)$ is diffeomorphic to 
$\big( V \times {\chi }^{-1}(s'), C^{\infty}({\mathfrak{A}}_{\mathrm{stg}}(V \times {\chi }^{-1}(s'))) \big)$. \medskip 
 
\noindent \textbf{Proof.} Our argument parallels the proof of lemma 8. The strong 
tangent map at $u \in U$ of the smooth mapping $\Psi $ (\ref{eq-twovnw}) is 
\begin{displaymath}
\begin{array}{l}
T_u\Psi : T_u U = T_vV \times T_y\big( {\chi }^{-1}(v) \big) \rightarrow 
T_vV \times T_{{\varphi }_v(y)}\big( {\chi }^{-1}(s') \big):  \\ 
\rule{0pt}{16pt}\hspace{1in} \mbox{\footnotesize $\begin{pmatrix} \xi \\ \eta \end{pmatrix}$} 
\mapsto \mbox{\footnotesize $\begin{pmatrix} I & D_1{\varphi }_v(y) \\ 0 & D_2{\varphi }_v \end{pmatrix} \,  
\begin{pmatrix} \xi \\ \eta \end{pmatrix}. $}
\end{array}
\end{displaymath}
For fixed $v \in V$ the mapping ${\varphi }_v$ is a diffeomorphism by hypothesis. So its 
partial derivative $D_2{\varphi }_v$ in the $T_y\big( {\chi}^{-1}(v) \big)$ direction is invertible. 
Thus $T_v\Psi $ is invertible, which implies that $\Psi $ (\ref{eq-twovnw}) is a smooth 
local diffeomorphism with strong tangent. $\Psi $ is a diffeomorphism because it is surjective and injective, using 
the same argument as in the proof of lemma 8.  \hfill $\square $ 

\section{Examples}

In this section we give some examples of a smooth mapping of subcartesian 
differential spaces, which is a smooth surjective submersion. 

\subsection{Preliminary example}

Let $E$ be a finite dimensional real vector space. Let $\Phi : H \times E \rightarrow E$ be a 
linear action of a compact subgroup $H$ of $\Gl (E, \R )$. The algebra $P(E)^H$ of $H$-invariant 
polynomials on $E$ is finitely generated, say by $p_1, \ldots , p_m$. Because the $H$ action 
$\Phi $ is linear we may choose $p_i$ to be homogeneous of degree $d_i >0$. We may 
assume that $m$ is minimal. Then ${\{ p_i \} }^m_{i=1}$ is a Hilbert basis of $P(E)^H$ and 
\begin{displaymath}
{\sigma }: E \rightarrow {\R }^m: x \mapsto \big( p_1(x), \ldots , p_m(x) \big) = \overline{x}
\end{displaymath}
is the Hilbert map corresponding to the action $\Phi $. Since elements of $P(E)^H$ separate 
orbits of the action $\Phi $ on $E$ and $\sigma $ is $H$-invariant, $\sigma $ induces a continuous 
bijective mapping 
\begin{displaymath}
\widetilde{\sigma }: E/H \rightarrow \Sigma = \sigma (E) \subseteq {\R }^m.
\end{displaymath}
$\Sigma $ is a semialgebraic subset of ${\R }^m$. The map $\widetilde{\sigma }$ is a 
homeomorphism from the locally compact Hausdorff orbit 
space $E/H$ onto $\Sigma $, see \cite[p.324]{cushman-bates}.  \medskip 

By Schwarz's theorem every smooth $H$-invariant function on $E$ is a smooth function of the 
$H$-invariant polynomials, that is for every $f \in C^{\infty}(E)^H$ there is an $F \in C^{\infty}({\R }^m)$ 
such that ${\sigma }^{\ast }F = f$. We say that $\overline{f} \in C^{\infty}(E/H)$ if and only if 
${\sigma }^{\ast}\overline{f} \in C^{\infty}(E)^H$. $C^{\infty}(E/H)$ is a differential structure on 
$E/H$ and the map
\begin{displaymath}
\sigma : (E, C^{\infty}(E)^H) \rightarrow (E/H, C^{\infty}(E/H) )
\end{displaymath}
is a smooth map of differential of spaces, since ${\sigma }^{\ast }\overline{f} \in C^{\infty}(E)^H$. 
The map $\sigma $ is surjective. The differential space $(E/H, C^{\infty}(E/H))$ is 
subcartesian, since the orbit space $E/H$ is a subset of ${\R }^m$. \medskip 

The strong tangent space ${\mathcal{T}}_xE$ at $x \in E$ is ${\spann }_{\R } \{ X(x) \} $, where 
$X \in \mathfrak{X}(E)^H$ is a smooth $H$-invariant vector field on $E$. The flow ${\varphi }_t$ of 
$X$ commutes with the $H$ action $\Phi $, that is, ${\Phi }_h \big( {\varphi }_t(x) \big) = 
{\varphi }_t\big( {\Phi }_h(x) \big)$ for every $(h,x) \in H \times E$ and every $t \in \R $. So 
${\varphi }_t$ induces a one parameter group ${\overline{\varphi }}_t$ of diffeomorphisms of 
$(E/H, C^{\infty}(E/H))$ into itself. ${\overline{\varphi }}_t$ is the flow of a vector field 
$\overline{X}$ on $E/H$. The vector field $\overline{X}$ is $\sigma $ related to the vector field $X$, because 
\begin{displaymath}
\overline{ {\varphi }_t(x)} = \sigma \big( {\varphi }_t(x) \big) ={\overline{\varphi }}_t\big( \sigma (x) \big) =
{\overline{\varphi }}_t(\overline{x}) . 
\end{displaymath}
Thus $\overline{X}(\overline{x}) \in {\mathcal{T}}_{\overline{x}}(E/H)$. Hence the strong tangent map 
\begin{displaymath}
{\mathcal{T}}_x\sigma : {\mathcal{T}}_xE \rightarrow {\mathcal{T}}_{\overline{x}}(E/H): 
X(x) \mapsto \overline{X}(\overline{x})
\end{displaymath}
of the Hilbert map is defined at every $x \in E$. ${\mathcal{T}}_x\sigma $ is 
surjective, because every smooth vector field on $E/H$ is $\sigma $ related to a smooth $H$-invariant 
vector field on $E$, see \cite{bates-cushman-sniatycki}. Since $\Sigma $ is a subset of 
${\R }^m$, the inclusion map $\iota : \Sigma \rightarrow {\R }^m$ has a strong tangent mapping. 
Hence $(E/H, {\mathfrak{A}}_{\mathrm{stg}}(E/H))$ is a strong subcartesian space. \medskip 

Because $E$ is a smooth manifold, the subcartesian differential space $(E, C^{\infty}(E)^H)$ has an atlas 
$\mathfrak{A} = \{ {\mathrm{id}}_E: E \rightarrow E\}$ consisting of one chart, which has a 
strong tangent. Thus $(E, {\mathfrak{A}}_{\mathrm{stg}}(E))$ is a strong subcartesian space. 
The Hilbert map $\sigma : (E, {\mathfrak{A}}_{\mathrm{stg}}(E)) \rightarrow 
(E/H, {\mathfrak{A}}_{\mathrm{stg}}(E/H))$ is a surjective smooth mapping with a strong tangent. 
Consequently, $\sigma $ is a surjective 
submersion of the differential space $(E, C^{\infty }(E)^H)$ onto the subcartesian differential space 
$(E/H, C^{\infty}(E/H))$. \medskip 

\subsection{Main example}

Let $M$ be a smooth manifold and let $G$ be a Lie group, which acts  properly on $M$ via the 
mapping $\Phi : G \times M \rightarrow M$. For each $x \in M$ let $G \cdot x = 
\{ {\Phi }_g(x) \in M \setrule \, g \in G \} $ be the orbit of the $G$ action through $x$. Let 
$M/G = \{ G \cdot x \setrule \, x \in M \}$ be the space of $G$ orbits on $M$ with $G$ orbit mapping
\begin{displaymath}
\pi : M \rightarrow M/G : x \mapsto G\cdot x = \overline{x}.
\end{displaymath}
The quotient topology on $M/G$ is Hausdorff. \medskip 

We say that a function $\overline{f}: M/G \rightarrow \R $ is smooth if and only if 
${\pi }^{\ast}\overline{f} \in C^{\infty}(M)^G$, the set of smooth $G$ invariant functions on $M$. 
Let $C^{\infty}(M/G)$ be the set of smooth functions on $M/G$. Then $C^{\infty}(M/G)$ is 
a differential structure on $M/G$, see \cite[p.321]{cushman-bates}. The differential space topology 
of $(M/G, C^{\infty}(M/G))$ is the same as the quotient topology, because smooth $G$ invariant 
functions on $M$ separate $G$ orbits of the action $\Phi $. Hence the differential space topology 
is Hausdorff. Since ${\pi }^{\ast }\big( C^{\infty}(M/G) \big) = C^{\infty}(M)^G$ by definition, the orbit map 
$\pi : (M, C^{\infty}(M)^G) \rightarrow (M/G, C^{\infty}(M/G))$ is a smooth surjective mapping of differential 
spaces. $\pi $ has a strong tangent map at $x \in M$, since its tangent map sends a $G$ invariant 
vector field $X$ on $M$ evaluated at $x$ to the vector field $\overline{X}$ on $M/G$ evaluated at 
$\pi (x) = \overline{x}$. Here $\overline{X}$ is induced by $X$. \medskip 

We now review the argument, which shows that $(M/G, C^{\infty}(M/G))$ is a subcartesian differential space. For 
each $x \in M$ let $G_x = \{ g \in G \setrule \, {\Phi }_g (x) = x \}$ be the isotropy group of 
the $G$ action $\Phi $ at $x$. For every $g \in G_x$ the linear transformation 
$T_x{\Phi }_g: T_xM \rightarrow T_xM$ leaves the subspace $T_x(G \cdot x)$ invariant. Hence we 
obtain the action 
\begin{displaymath}
\ast : H \times E \rightarrow E: (h,v) \mapsto  T_x{\Phi }_hv 
\end{displaymath}
of a compact group $H = G_x$ on a real finite dimensional vector space $E= T_xM/T_x(G\cdot x)$. 
Let $B$ be an $H$-invariant open neighborhood of $0$ in $E$. On $G \times B$ we have an action of $H$ 
given by 
\begin{displaymath}
\mu : H \times (G \times B) \rightarrow G \times B: \big( h, (g,b) \big) \mapsto (gh^{-1}, h \ast b), 
\end{displaymath}
which is free and proper. Therefore the $H$ orbit space $G{\times }_H B$ is a smooth manifold. The 
$G$ action on $G \times B$ defined by 
\begin{displaymath}
\nu : G \times (G \times B) \rightarrow G \times B: \big( g', (g,b) \big) \mapsto (g'g,b)
\end{displaymath}
commutes with the $H$ action $\mu $. Hence it induces a $G$ action on $G{\times }_HB$. Thus 
there is a $G$ invariant open neighborhood $U$ of $x$ in $M$, an $H$-invariant open neighborhood 
$B$ of $0$ in $E$ and a diffeomorphism $\varphi : G{\times }_HB \rightarrow U$, which intertwines the $G$ action on $G{\times }_HB$ with the $G$ action on $U$, see \cite[p.318]{cushman-bates}. In fact, we have 
a diffeomorphism $\phi $ sending the differential space $(B/H, C^{\infty}(B/H))$ onto the 
differential space $(U/G , C^{\infty}(U/G))$, see \cite[p.323]{cushman-bates}. Consider the 
diagram \bigskip 

\hspace{1.6in}\begin{tikzcd}
B \arrow{r}{\psi} \arrow{d}{\rho |B} & U \arrow{d}{\pi |U} \\
B/H \arrow{r}{\phi } & U/G
\end{tikzcd} 
\vspace{.2in}
\par \noindent Here $\psi $ is the diffeomorphism given by Bochner's lemma \cite[p.306]{cushman-bates}, 
$\rho |B$ is the orbit map of the $H$ action $\ast$ on $E$ restricted to $B$, and 
$\pi |U$ is the orbit map $\pi $ of the $G$ action on $M$ restricted to $U$. Since the smooth 
maps $\rho |B$, $\psi $, and $\pi |U$ have a strong tangent and the diagram commutes, it 
follows that the diffeomorphism $\phi $ has a strong tangent. \medskip 

We now argue that the orbit space of a proper action of $G$ on a smooth manifold $M$ is 
a subcartesian differential space. Let $x \in M$. Consider the action $\ast $ of $H = G_x$ on the 
finite dimensional real vector space $E = T_xM/T_x(G\cdot x)$. Let 
\begin{displaymath}
\sigma : E \rightarrow {\R }^m: x \mapsto \big( p_1(x), \ldots , p_m(x) \big)
\end{displaymath}
be the Hilbert map corresponding to the Hilbert basis ${\{ p_i \} }^m_{i=1}$ of the space 
$P(E)^H$ of $H$-invariant polynomials on $E$. Since $\sigma $ is $H$-invariant, it induces the map 
\begin{displaymath}
\widetilde{\sigma }: E/H \rightarrow \Sigma = \sigma (E) \subseteq {\R }^m, 
\end{displaymath}
which is a homeomorphism and a smooth surjective diffeomorphism of the differential spaces 
$(E/H, C^{\infty}(E/H))$ and $(\sigma (E), C^{\infty}(\sigma (E)))$ $ \subseteq 
({\R }^m, C^{\infty}({\R }^m))$ with a strong tangent. \medskip 

Let $\mathcal{F}$ be the family of functions $\overline{f}: \Sigma \subseteq 
{\R }^m \rightarrow \R $ such that 
${\sigma }^{\ast }\overline{f}: E \rightarrow \R$ is a smooth function on $E$. Note that 
${\sigma }^{\ast }\overline{f}$ is $H$-invariant. Let $C^{\infty}(\Sigma )$ be the differential structure on 
$\Sigma $ generated by $\mathcal{F}$. Let $\widetilde{\mathcal{F}}$ be the family of function $\widetilde{f}: \Sigma \subseteq {\R }^m \rightarrow \R$ such that there is an $F \in C^{\infty}({\R }^m)$ with $F|\Sigma = \widetilde{f}$. Let ${\widetilde{C}}^{\infty}(\Sigma )$ be the differential structure on $\Sigma $ generated by 
$\widetilde{\mathcal{F}}$. The differential spaces $(\Sigma , C^{\infty}(\Sigma ))$ and 
$(\Sigma , {\widetilde{C}}^{\infty}(\Sigma ))$ are diffeomorphic. The differential space topology on $\Sigma $ generated by $C^{\infty}(\Sigma )$ is the same as the topology on $\Sigma $ induced from ${\R }^m$.
\medskip 

We now show that $(\Sigma , C^{\infty}(\Sigma )) \subseteq ({\R }^m, C^{\infty}({\R }^m))$ is a 
locally compact subcartesian differential space. Shrinking $B$ and $U$ if necessary, we may 
assume that $B = \{ x \in E \setrule \, \beta (x,x) < c \}$ for some $c>0$. Here $\beta $ 
is an $H$-invariant inner product on $E$. The map 
$\widetilde{\sigma }$ is a homeomorphism of the orbit space $B/H$ onto 
$\sigma (B) = \{  y \in \Sigma = \sigma (E) \setrule \, P(y) < c \}$, where 
$P$ is a polynomial on ${\R }^m$ such that $\beta (x,x) = P\big( \sigma (x) \big)$ for 
all $x \in E$. $\sigma (B)$ is an open 
semialgebraic subset of the closed semialgebraic set $\Sigma $. Let \linebreak 
$\phi : (B/H, C^{\infty}(B/H)) \rightarrow (U/G , C^{\infty}(U/G))$ be the diffeomorphism in the above 
diagram. The map 
\begin{displaymath}
\widetilde{\sigma }\comp {\phi }^{-1}: (U/G , C^{\infty}(U/G)) \rightarrow \big( \sigma (B), 
C^{\infty}(\sigma (B)) \big) \subseteq ({\R }^m, C^{\infty}({\R }^m) ) 
\end{displaymath}
is a diffeomorphism, see \cite[p.327]{cushman-bates}. Because the open sets $\{ U/G \}$ cover the orbit space $M/G$, it follows that $(M/G, C^{\infty}(M/G))$ is a subcartesian differential space. 
Since the mappings $\widetilde{\sigma }$ and $\phi $ are smooth diffeomorphisms with a strong tangent, it follows that the chart map $\widetilde{\sigma } \comp {\phi }^{-1}: U/G \rightarrow B/H \subseteq {\R }^m$ is a 
smooth diffeomorphism with a strong tangent. Hence $(M/G, {\mathfrak{A}}_{\mathrm{stg}}(M/G))$ is a strong subcartesian space.

\end{document}